\documentclass[12pt]{article}
\usepackage{amsfonts}
\usepackage{amssymb}
\usepackage{mathrsfs}
\usepackage{float}
\usepackage{fancyhdr}
\usepackage{amsmath,,amsthm}
\usepackage{amsfonts,amssymb}
\usepackage{color}
\usepackage{booktabs,tabularx}
\usepackage{dsfont}
\usepackage{epsfig}
\usepackage{graphicx}
\usepackage{subfigure}
\renewcommand{\thesubfigure}{\roman{subfigure}}
\makeatletter
\renewcommand{\@thesubfigure}{(\thesubfigure)\space}
\renewcommand{\p@subfigure}{\thefigure}
\makeatother

\date{ }

\title{Convergence of generalized Collatz problem in $k$-adic field}
\author{Yushu Zhu, Sensen Chen and Qing-You Sun\footnote{Corresponding author:
qysun@hznu.edu.cn}
\\ \small Hangzhou Normal University, Hangzhou
311121, China }

\begin{document}
\maketitle
\begin{abstract}
In this article, we define a new $k$-adic series transformation
called $\mathcal{Z}$-transformation and probe into its fixed point
and periodicity. We extend the number field of the transform period
problem to a wider $k$-adic field. Different constraints are imposed
on $k$, then different periodic columns are formed after finite
$\mathcal{Z}$ transformations. We obtain that their periodic
sequences are $M_1= \{1,2\}$ and $M_2=\{1,2\}\cup \{n_0\} \cup
\{n'\}$ respectively after derivation. As an application, it can
provide a reference for C problems in more complex algebraic
systems.

\vskip 3mm



\end{abstract}

\baselineskip=5mm

\section{Introduction}

\hspace{0.5cm} The Collatz conjecture is named after Lothar Collatz,
who introduced the idea in 1937. It has been widely studied in the
past 100 years, and many achievements with great value have been
obtained, although the Collatz problem cannot be effectively solved.

The results of the Collatz problem can be applied to modern
cryptography, and can be easily transformed into graph theory, and
then extended to a wider range of algebraic and algebraic geometry.
The essence of the Collatz problem is a fixed-point problem, its
research methods and results are of great significance to the study
of fixed point problems, and thus have an important impact on the
development of modern mathematics. The problem also has important
applications in power systems and fractal geometries.

The Collatz problem is concerned with the function $f:\mathds{Z}^+
\rightarrow \mathds{Z}^+$ defined as \cite{1}

\begin{equation} \label{eq:1}
f(n)=\left\{\begin{array}{c}{\displaystyle\frac{n}{2}, \text { when
} n \text { is even, }}
\\ {\displaystyle\frac{3 n + 1}{2}, \text { when } n \text { is odd.
}}\end{array}\right.
\end{equation}

The famous Collatz problem asserts that iteration series $f(n),f^2
(n),f^3 (n),\cdots$ of every positive number $n$ will eventually
reach the integer one. This conjecture is still an unsolved mystery,
many articles have studied this issue.

Matthews and Watts \cite{2} gave generalized Collatz maps

\begin{equation*}
f_{d,m,r}(x)=\displaystyle\frac{m_i x +r_i}{d},\text { if } x\equiv
i(\bmod d),
\end{equation*}
where $d>1$ is an integer, $m$ and $r$ are $d$-dimensional vectors
such that $m_i x+r_i\equiv 0(\bmod d)$ for all $i$. The number $n$
is supposed to be cyclic for $f_{d,m,r}$ if the series
$\left\{f^i_{d,m,r}(n)\right\}_{i\in \mathds{N}}$ is a periodic
column. As is known, we can see that the density of the set of the
integer will stop for $1$ after finite times of Collatz maps which
lead us to seek for other periodic column.

Certainly, these problems and their promotion problems have many
research results, especially through computers to solve these
problems (\cite{3}-\cite{6}). The results of this article are mainly
inspired by \cite{7} and \cite{8}, which extend the problem to
$p$-adic. This led us to consider this problem from a deeper
perspective of domain theory, which makes us want to study a class
of number theory problem.

In section 2 and 3, we introduce some notations and definitions, and
give our main results. In section 4, by studying the properties of
the $\mathcal{Z}$-transform, Lemma 2 mentions that when the digit of
$n_0$ is greater than $2$, then after finite
$\mathcal{Z}$-transformations, it can degenerate into a single-digit
number or a double-digit number in $k$-adic. So, we mainly discuss
these two cases in the proof of Theorem 1. For both cases, we
discuss the result of $n_0$ after finite
$\mathcal{Z}$-transformations. Using the mathematical induction
method, we obtain that after finite $\mathcal{Z}$-transformations
under the hypothesis, there is only one periodic column
$M_1=\{1,2\}$. In Theorem 2, we weaken the restriction on $k$, then
proved that the period column of the $\mathcal{Z}$-transformation is
listed as $M_2=\{1,2\}\cup\{n_0 \}\cup\{n'\}$.

\section{Notation and Definition}

\hspace{0.5cm} For any given $m$-digit integer expressed in $k$-adic

\begin{equation} \label{eq:n}
n_0=a_{m-1} k^{m-1}+a_{m-2} k^{m-2}+\dots+a_{1} k+a_{0},
\end{equation}
where $a_i \in \mathds{Z}_0^+$, $0 \leq a_i < k$ $(0 \leq i \leq
m-1)$, and $a_{m-1} \neq 0$.

We construct the function in $k$-adic,
\begin{equation} \label{eq:f}
f(n)=\left\{\begin{array}{cc}{\displaystyle\frac{(n+p-1)(n+2 p-1)}{p
\cdot p}} & {n \equiv 1(\bmod p)},
\\ {\displaystyle\frac{n+p-2}{p}} & {n \equiv 2(\bmod p)},
\\\vdots
\\{\displaystyle\frac{n+1}{p}} & {n = p-1(\bmod p)},
\\ {\displaystyle\frac{n}{p}} & {n \equiv
0(\bmod p)}\end{array}\right.
\end{equation}
where $0\leq n <k$.

\noindent\textbf{Definition 1}\quad \emph{ Let
$\mathcal{Z}(n_0)=\displaystyle\sum_{i=0}^{m-1}f(a_i)$, where $n_0$,
$a_i$ and $f$ are defined as above, we call it
$\mathcal{Z}$-transformation of $n_0$ in $k$-adic.}

\noindent\textbf{Definition 2}\quad \emph{ Transformation sequence
$\{\mathcal{Z}^{i}(n_0)\}, \quad (i\in \mathds{Z}_0^+)$, is called
$\mathcal{Z}$-transformation sequence in $k$-adic.}

Denote $n_i=\mathcal{Z}^{i}(n_0)$, in which
$\mathcal{Z}^{0}(n_0)=n_0$, that is, $\{n_{i}\}, \quad (i\in
\mathds{Z}_0^+)$, is also the $\mathcal{Z}$-transformation sequence.
It's clear that $n_i\in \mathds{N}$.

Our article mainly discusses the periodic column of
$\mathcal{Z}$-transform sequence.

\section{The Main Results}

\noindent\textbf{\emph Hypothesis:}\quad \emph{\\
(a) $2p-1 \leq k \leq 3p^2$, when $p\geq 3$;\\
(b) For any nonnegative integer $q$ less than $\sqrt{k}$, it
satisfies
$(q+1)(q+2) \not\equiv -1(\bmod p)$.\\
(c) For any nonnegative integer $q$ less than $\sqrt{k}$, it
satisfies $(q+1)(q+2) \not\equiv k(\bmod p)$, and $k\neq
(q+1)(q+2)-1-qp$.}

Under the Hypothesis, we give our main result as the following two
Theorem.

\noindent\textbf{Theorem 1}\quad \emph{ If the Hypothesis is
established, then after finite $\mathcal{Z}$-transformations,
$\mathcal{Z}$-transformed sequence $\{n_i \}, i\in \mathds{N}$ only
has one period column $M_1 = \{1,2 \}$.}

\noindent\textbf{Theorem 2}\quad \emph{ If only (a) in the
Hypothesis is hold, then after finite $\mathcal{Z}$-transformations,
$\mathcal{Z}$-transformed sequence $\{n_i \}, i\in \mathds{N}$ has
the period column $M_2=\{1,2\}\cup \{n_0\} \cup \{n'\}$, where $n'$
is a integer which satisfies $\mathcal{Z}(n')=n'$.}

\section{The Proof of the Theorems}

\hspace{0.5cm} To prove the Theorems, we give two Lemmas first.

\noindent\textbf{Lemma 1}\quad \emph{The result of adding two
$k$-adic numbers is the same as the result of adding them in decimal
to become $k$-adic.}

This conclusion is obvious (see \cite{2}). We omit it.

\noindent\textbf{Lemma 2}\quad \emph{For any given $m$-digit
positive integer $n_0$ in $k$-adic, under the Hypothesis,
$\mathcal{Z}(n_0)$ is a integer not exceeding $(m-1)$-dgit in
$k$-adic when $m\geq 3$.}

\begin{proof}
Denote
\begin{equation*}
k = tp +s +1,
\end{equation*}
where $r,s\in \mathds{N}$ and $1\leq s\leq p$.

Due to $0\leq a_i \leq k-1$, by the definition of $f(n)$ in
\eqref{eq:f}, it obtains
\begin{equation*}
f(a_i)\leq f(tp+1)=(t+1)(t+2).
\end{equation*}
Thus, $\mathcal{Z}(n_0) \leq m(t+1)(t+2)$.

For $m=3$, by (a) in the Hypothesis, we have
\begin{equation*}
\mathcal{Z}(n_0) \leq 3(t+1)(t+2) < k^{3-1} .
\end{equation*}
For $m=j$, we assume that it holds
\begin{equation*}
\mathcal{Z}(n_0) \leq j(t+1)(t+2) < k^{j-1} .
\end{equation*}
Then, for $m=j+1$, it's obviously that
\begin{equation*}
\mathcal{Z}(n_0) \leq (j+1)(t+1)(t+2) < k^{j-1}\cdot
\displaystyle\frac{j+1}{j} < k^{(j+1)-1} .
\end{equation*}
Therefore, by mathematical induction, we have
\begin{equation}
\mathcal{Z}(n_0) < k^{m-1} ,
\end{equation}
where $m \geq 3$. That means, $\mathcal{Z}(n_0)$ is a integer not
exceeding $(m-1)$-dgit in $k$-adic when $m\geq 3$. Thus, the proof
of Lemma 2 is done.

\end{proof}

\noindent{\bf \emph{Proof of Theorem 1.}} For any given $m$-digit
positive integer $n_0$ in $k$-adic, we divide the proof of Theorem 1
into three parts. Here, we can omit the case $n_0 = 1$, which make
Theorem 1 be true obviously.

{\bf Case 1.} $n_0$ is a single-digit number.

In this case, it is easy to show that $1 < n_0 < k$, and
$\mathcal{Z}(n_0) = f(n_0)$.

{\bf(I).} If $n_0 \not\equiv 1(\bmod p)$, by the definition of
$f(n)$, we can easily get $\mathcal{Z}(n_0) = f(n_0) < n_0$. It
implies that the $\mathcal{Z}$ transformation make
$\mathcal{Z}(n_0)$ smaller than $n_0$, and $\mathcal{Z}(n_0)$ is
still a single-digit number. Therefore, there exists a nonnegative
integer $\lambda$, which satisfies $\mathcal{Z}^{\lambda}(n_0) = 1$,
or $\mathcal{Z}^{\lambda}(n_0) \equiv 1(\bmod p)$. The second case
is just what we will discuss next.

Analogously, in the rest of the proof, we just need to prove that
there exists a finite integer $\lambda$ make
$\mathcal{Z}^{\lambda}(n_0) < n_0$. Since we can think
$\mathcal{Z}^{\lambda}(n_0)$ as a new $n_0$, and do the same
discussion again and again, until $\mathcal{Z}^{\lambda}(n_0)=1$.

{\bf(II).} If $n_0 \equiv 1(\bmod p)$ and $n_0 \neq 1$, let
$n_0=rp+1$. So, it can obtain that
\begin{equation*}
\mathcal{Z}(n_0) = f(n_0) = (r+1)(r+2).
\end{equation*}

Assume that $\mathcal{Z}(n_0)=a_1 k+b_1$, where $a_1,b_1 \in
\mathds{Z}_0^+$ and $b_1<k$. Using (a) in the Hypothesis and the
definition of $f$, it can be obtained by direct calculation that
$a_1 \leq 3$.

{\bf(i).} When $a_1=0$, that is $\mathcal{Z}(n_0) =
f(n_0)=(r+1)(r+2)<k$.

If $(r+1)(r+2) \not\equiv 1(\bmod p)$, we deduce that
\begin{equation*}
\mathcal{Z}^2(n_0) =
\mathcal{Z}((r+1)(r+2))=\frac{(r+1)(r+2)+\theta}{p},
\end{equation*}
where $\theta$ is a nonnegative integer less than $p$. Noticing
$n_0>1$ which we discussing here, then by direct calculation,
$\mathcal{Z}^2(n_0) \geq n_0$ is only available when $n_0>p(p^2-3)$.
On another hand, $f(n_0) < k$ implies $n_0<p(p-3)+1<p(p^2-3)$ when
$p>7$. Thus, it indicates that $\mathcal{Z}^2(n_0) < n_0$.

If $(r+1)(r+2) \equiv 1(\bmod p)$, and $\mathcal{Z}(n_0)> n_0$ are
established, then for the finite number $k$, there exists $\lambda$
such that $\mathcal{Z}^\lambda(n_0)> k$, we will discuss it in (ii)
and (iii). What we want to explain here is that
$\mathcal{Z}^\lambda(n_0)>\mathcal{Z}^{\lambda-1}(n_0)$ in this
case. Also, for the integer $n_0$, $\mathcal{Z}(n_0)=n_0$ does not
hold.

{\bf(ii).} When $a_1=1$, that is $\mathcal{Z}(n_0) = k+b_1$.

If $b_1 \not\equiv 1(\bmod p)$, we have
\begin{equation*}
\mathcal{Z}^2(n_0) =
f(1)+f(b_1)=2+\frac{b_1+\theta}{p}<\frac{k}{p}+3,
\end{equation*}
where $\theta$ is a nonnegative integer less than $p$. Noting (a) in
the Hypothesis, it follows that
\begin{equation*}
\mathcal{Z}^2(n_0)< \frac{k}{p}+3< k\leq n_0.
\end{equation*}

If $b_1 \equiv 1(\bmod p)$, then
\begin{equation*}
\mathcal{Z}^2(n_0) = f(1)+f(b_1)=2+f(b_1)\triangleq a_2 k+b_2.
\end{equation*}
For $\mathcal{Z}^i(n_0),i>2$, we can make the similar notation. And
we only discuss the case $b_i \equiv 1(\bmod p)$, otherwise, it will
be the case we discussed above. Denote $b_i = q_i p +1$. $b_i<k$
means $q_i<\frac{k}{p}$. Similar to the discussion about $a_1$, it
can be obtained by direct calculation that $a_2 \leq 3$.

{\bf(A).} Suppose that $\mathcal{Z}^2(n_0)=b_2<k$ is established, we
only need to deal with the case $\mathcal{Z}^2(n_0)\geq n_0$. If
$\mathcal{Z}^2(n_0)> n_0$, then by the definition of $f(n)$, for the
finite number $k$, there exists $\lambda$ such that
$\mathcal{Z}^\lambda(n_0)> k$, we will discuss it in the following.
On the other hand, we have $\mathcal{Z}^2(n_0)\neq n_0$ when $q_1$
satisfies (b) in the Hypothesis, since
$\mathcal{Z}^2(n_0)=f(b_1)+2=(q_1+1)(q_1+2)+2\neq rp+1$.

{\bf(B).} Suppose that $\mathcal{Z}^2(n_0)=k+b_2$ is established,
then we have $b_1 \neq b_2$. Otherwise,
$f(b_1)+2=\mathcal{Z}^2(n_0)=\mathcal{Z}(n_0)=f(n_0)$, that is
$(q_1+1)(q_1+2)+2=(r+1)(r+2)$. This can not be established for any
nonnegative integers $r$ and $q_1$.

If $b_2<b_1$, then $f(b_2)<f(b_1)$, that will make
$\mathcal{Z}^i(n_0) < \mathcal{Z}^{i-1}(n_0)$, until
$\mathcal{Z}^i(n_0)<k$. Also, if $q_{i-1}$ satisfies (b) in the
Hypothesis, it's easy to get $\mathcal{Z}^i(n_0) \neq n_0$. Then,
there are finite $i$ which satisfy $n_0 < \mathcal{Z}^i(n_0)< k$ and
different from each other. Thus, there exists $\lambda$ such that
$\mathcal{Z}^\lambda(n_0)< n_0$ or $b_{\lambda+1} > b_\lambda$.

If $b_2>b_1$, that will make $\mathcal{Z}^i(n_0) >
\mathcal{Z}^{i-1}(n_0)$, until $\mathcal{Z}^i(n_0)\geq 2k$. We will
discuss it next.

{\bf(C).} Suppose that $\mathcal{Z}^2(n_0)=2k+b_2$ is established.

If $b_2 < b_1$, then
$\mathcal{Z}^3(n_0)=1+f(b_2)<2+f(b_1)=\mathcal{Z}^2(n_0)$. Thus,
there exists $i \geq 2$ such that $\mathcal{Z}^i(n_0)=k+b_i$ or
$\mathcal{Z}^i(n_0)<k$, which we have discussed in (A) and (B). The
only difference is the way to prove $\mathcal{Z}^i(n_0)\neq n_0$.
Assume $\mathcal{Z}^i(n_0)=f(2)+f(b_{i-1})=1+(q_{i-1}+1)(q_{i-1}+2)=
n_0$, we deduce $(q_{i-1}+1)(q_{i-1}+2)\equiv 0(\bmod p)$ since
$n_0\equiv 1(\bmod p)$. It means $q_{i-1}+1=p$, or $q_{i-1}+2=p$
when $q_{i-1}<\sqrt{k}$. It implies $n_0=p(p+1)+1$, or
$n_0=(p-1)p+1$. That is, $\mathcal{Z}(n_0)=(p+2)(p+3)$, or
$\mathcal{Z}(n_0)=p(p+1)$. Noting $n_0<k<\mathcal{Z}(n_0)$, from the
first one we have $\mathcal{Z}(n_0)=\left(p(p+1)+2\right)+(4p+4)$,
then $\mathcal{Z}^2(n_0)\leq f(1)+f(4p+1)=21\leq
2\left(p(p+1)+2\right)\leq 2k$ when $p\geq 3$. And from the seconde
one we have $\mathcal{Z}(n_0)=\left(p(p-1)+2\right)+(2p-2)$, then
$\mathcal{Z}^2(n_0)\leq f(1)+f(p+1)=8\leq p(p-1)+2\leq k$ when
$p\geq 3$. These two case both contradict
$\mathcal{Z}^2(n_0)=2k+b_2>2k$.

If $b_2 = b_1$, then
$\left(2+f(b_1)\right)-\left(1+f(b_2)\right)=1$, so
$\mathcal{Z}^3(n_0)\neq a_3k +b_3$ when $\mathcal{Z}^2(n_0)=a_2k
+b_2$, where $b_i = q_i p+1$. That will be the case $b_i \not\equiv
1(\bmod p)$.

If $b_2>b_1$, that will make $\mathcal{Z}^i(n_0) >
\mathcal{Z}^{i-1}(n_0)$, until $\mathcal{Z}^i(n_0)\geq 3k$. We will
discuss it next.

{\bf(D).} Suppose that $\mathcal{Z}^2(n_0)=3k+b_2$ is established,
this will be similar to $a_1=3$. We will discuss it in (iv).

{\bf(iii).} When $a_1=2$, that is $\mathcal{Z}(n_0) = 2k+b_1$, which
deduces $\mathcal{Z}^2(n_0)=1+f(b_1)$. Then we can make a similar
discussion as what we do in (ii).

It should be noted that, in this case we need the condition
$(q_i+1)(q_i+2) \not\equiv k(\bmod p)$, that is the first condition
of (c) in the Hypothesis, to obtain $b_1\neq b_2$ corresponding to
what we did in (B) of (ii).

{\bf(iv).} When $a_1=3$, that is $\mathcal{Z}(n_0) = 3k+b_1$.

Since $n_0=rp+1$ and $n_0<k$, then by solving
$\mathcal{Z}(n_0)=(r+1)(r+2) \geq 3k \geq 3rp+6$, it implies that
$n_0=3p^2-p+1$ or $n_0=3p^2-2p+1$.

If $n_0=3p^2-p+1$, then $\mathcal{Z}(n_0)=f(n_0)=3p(3p+1)=3k+b_1$,
and $k \geq 3p^2-p+2$. This indicates $b_1\leq 6p-6$. Noting $p\geq
3$, when $p\geq 4$, $\mathcal{Z}^2(n_0)=f(3)+f(b_1)\leq 43<n_0$. And
when $p=3$, $\mathcal{Z}^2(n_0)=f(3)+f(b_1)\leq 23<n_0$.

If $n_0=3p^2-p+1$, we can get the similar conclusion by the same
direct calculation.

In summary, we have a brief proof of Theorem 1 for Case 1.

{\bf Case 2.} $n$ is a 2-digit number.

Denote $n_0=c_0 k+d_0$, then
$\mathcal{Z}_(n_0)=f(c_0)+f(d_0)\triangleq c_1 k+d_1$. By what we
discuss in (II) of Case 1, since $c_0,d_0<k$, it will imply
$f(c_0),f(d_0)\leq 3k + (6p-6)$. Thus, $\mathcal{Z}_(n_0)=c_1
k+d_1<7k+d_1$, then, $\mathcal{Z}^2(n_0)=f(c_1)+f(d_1)\leq
2+f(d_1)<3k+(6p-4)$.

If $d_1\not\equiv 1(\bmod p)$, the conclusion is clearly right.

If $d_1\equiv 1(\bmod p)$, think $2+f(d_1)$ as $a_1k+b_1$ in (II) of
Case 1, then by the same way we did, there exists a nonnegative
integer $\lambda$, which satisfies
$\mathcal{Z}^{\lambda}\left(\mathcal{Z}(n_0)\right) = 1$. In this
case, we will use (c) in the Hypothesis to make sure that
$\mathcal{Z}^{\lambda}\left(\mathcal{Z}(n_0)\right) \neq n_0$.

Thus, we finish the proof of Theorem 1 for Case 2.

{\bf Case 3.} $n$ is a m-digit number, where $m \geq 3$.

In this case, by Lemma 2, we can easily get that there exists a
nonnegative integer $\lambda$, $\lambda\leq m-2$, which satisfies
$\mathcal{Z}_k^{\lambda}(n)$ is a 2-digit number. It will become
Case 2.

Overall, it is not difficult to see that there will exist a
nonnegative integer $\lambda$, which satisfies
$\mathcal{Z}^{\lambda}(n) =1$. And noting $\mathcal{Z}(1)=f(1)=2$,
$\mathcal{Z}(2)=f(2)=1$, therefore, for any positive integer $\mu$,
which satisfies $\mu \geq \lambda$, we have $\mathcal{Z}^{\mu}(n)
\in \{1, 2\}$.

Hence, Theorem 1 is
proved.\\
\rightline{$\Box$}

\noindent{\bf \emph{Proof of Theorem 2.}} The proof of Theorem 2 is
similar to the proof of Theorem 1. The only difference is that,
without (c) in the Hypothesis, we can not obtain
$\mathcal{Z}(n_0)\neq n_0$ and $\mathcal{Z}^{i}(b_i)\neq
\mathcal{Z}^{i-1}(b_i)$. Therefor, we need to add $ \{n_0\} \cup
\{n'\}$ to the period column, where $n'$ is a integer which
satisfies $\mathcal{Z}(n')=n'$.

Hence, Theorem 2 is proved.\\
\rightline{$\Box$}

\section{Examples}

\hspace{0.5cm} In this section, we give some examples for different
$n$, $k$, and $p$, to see the periodic characteristic of the
$\mathcal{Z}$ transformation.

\noindent{\bf \emph{Example 1.}} Take $n=123789$, $k=137$, and
$p=11$, this set satisfies the requirement of the Hypothesis, that
is Theorem 1.

We can see from Figure 1 that by three times of the $\mathcal{Z}$
transformation, $n=123789$ will become $1$ in $137$-adic.

\noindent{\bf \emph{Example 2.}} Take $n=9827$, $k=5$, and $p=3$,
this set satisfies (a) in the Hypothesis. But when $q=0<\sqrt{5}$,
it doesn't satisfies (b).

\noindent{\bf \emph{Example 3.}} Take $n=8512$, $k=10$, and $p=5$,
this set also satisfies (a) in the Hypothesis. But when
$q=3<\sqrt{10}$, the first requirement in (c) is not satisfied.

Example 2 and 3 are still satisfies the requirement of Theorem 2.
From Figure 2, we can see that $\{n'\}$ in Theorem 2 for them are
$\{4,6\}$ and $\{6\}$ respectively.


\begin{figure}[htbp]
\centering
\includegraphics[width=3.1in]{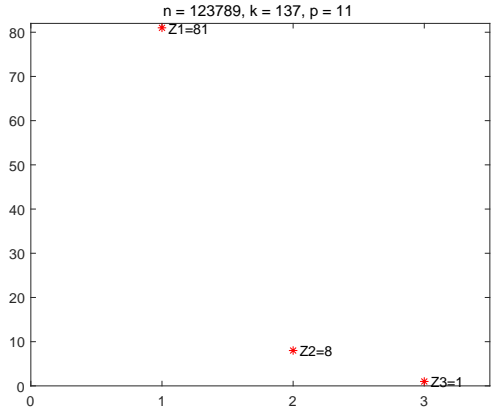}
\centering \caption{A set of $n, k, p$, which satisfies the
Hypothesis. }
\end{figure}

\begin{figure}[htbp]
\centering
\subfigure[Doesn't satisfy (b) in the Hypothesis.]{
\begin{minipage}[t]{0.45\linewidth}
\centering
\includegraphics[width=3.1in]{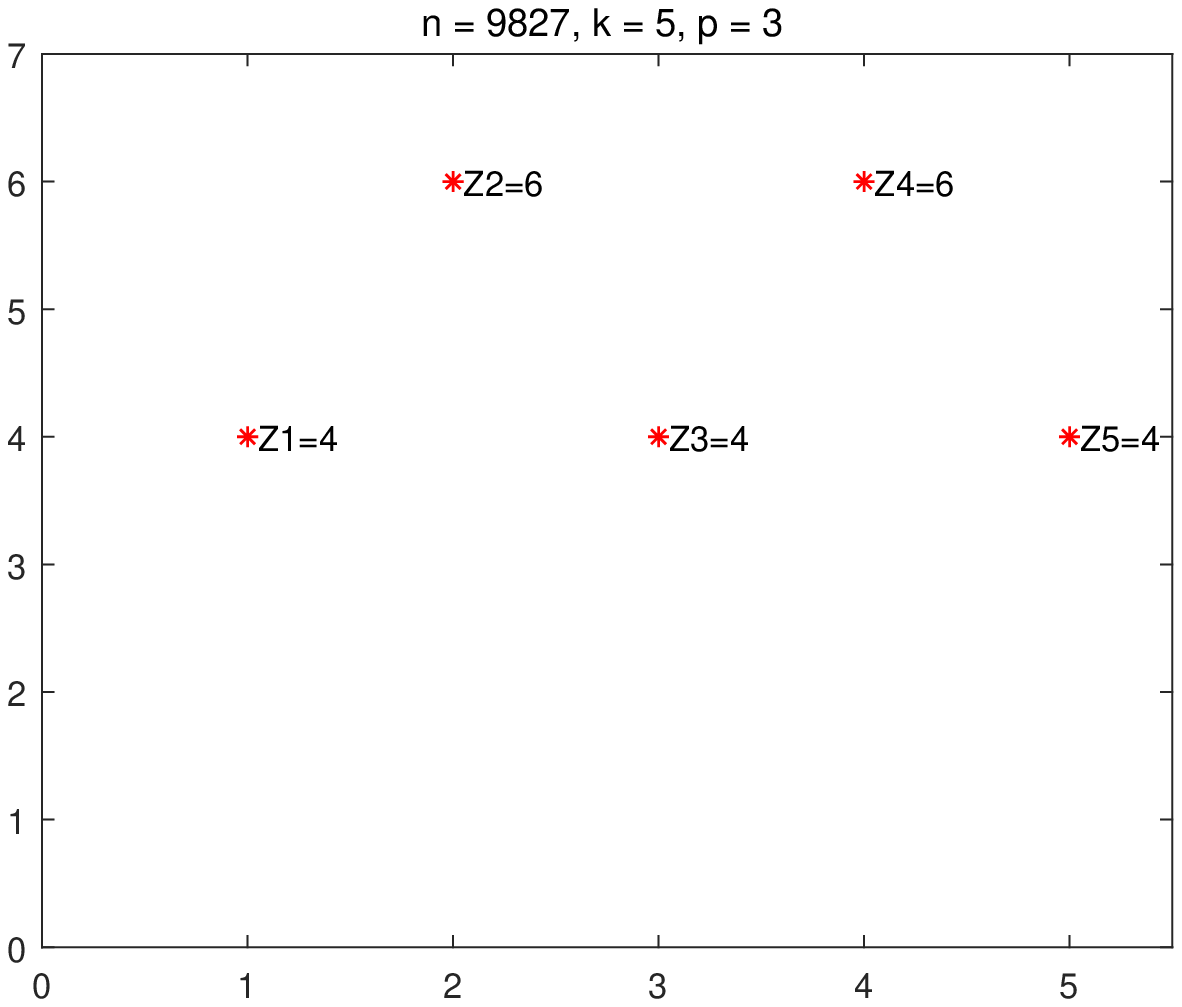}
\end{minipage}%
}%
\hspace{0.8cm}
\subfigure[Doesn't satisfy (c) in the Hypothesis.]{
\begin{minipage}[t]{0.45\linewidth}
\centering
\includegraphics[width=3.1in]{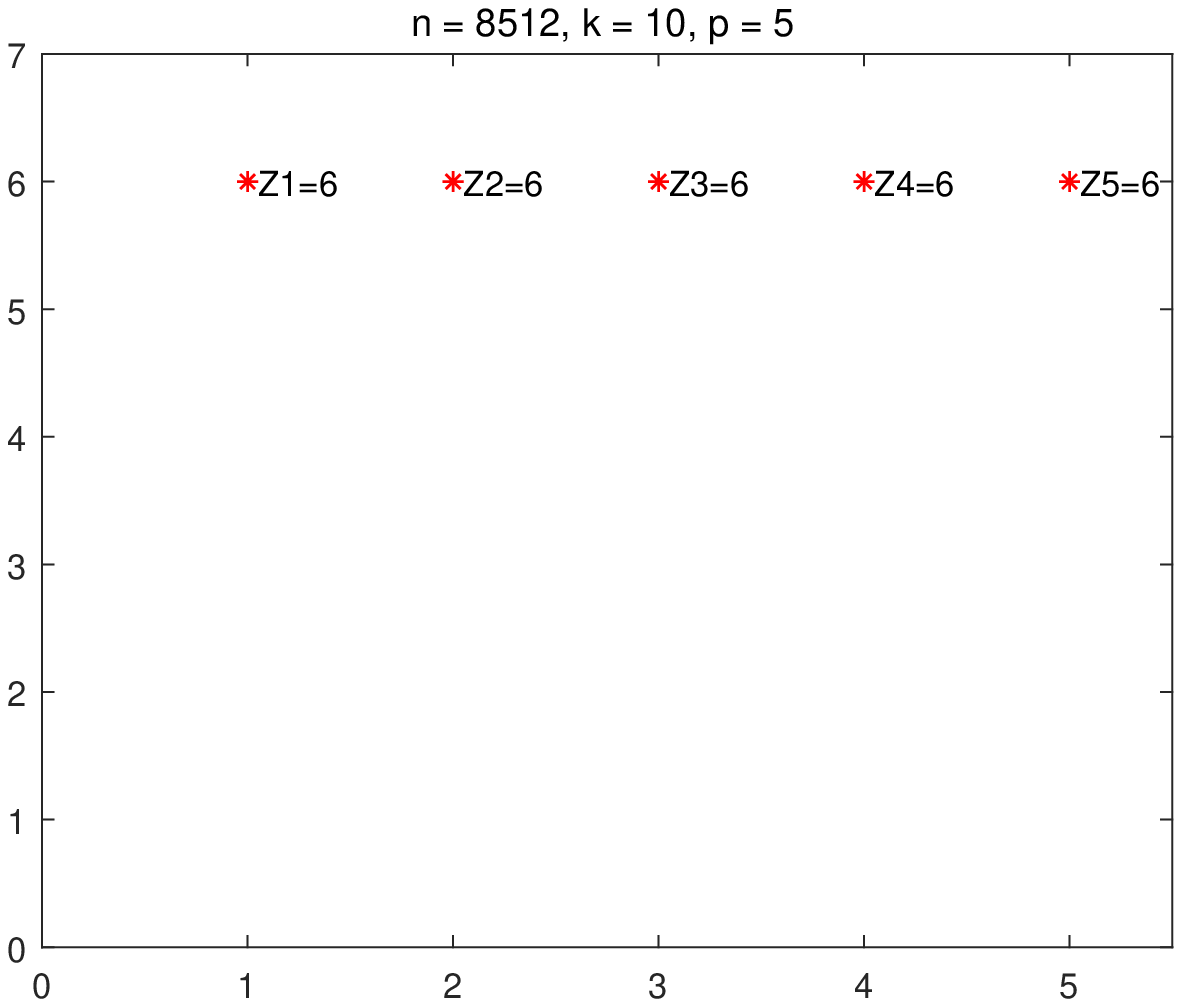}
\end{minipage}
}%
\centering \caption{Two sets of $n, k, p$, which satisfy the
requirement of Theorem 2.}
\end{figure}
\section*{Acknowledgements}

\hspace{0.5cm} This article is supported by Top Disciplines(Class-A)
of Zhejiang Province and Teaching Reform Project of Hangzhou Normal
University.


\end{document}